\documentclass[12pt]{amsart}
\usepackage{amsmath}
\newtheorem{theorem}{Theorem}
\newtheorem{lemma}{Lemma}

\numberwithin{equation}{section}

\begin{document}
\title{A converse theorem for $\Gamma_0(13)$}
\author{J.B. Conrey, David W. Farmer, B.E. Odgers,  and N.C. Snaith}
\thanks{The initial ideas for this work came about at the workshop ``Converse Theorems"
held at Bristol University in April, 2005.
The authors thank Bristol for making that work-week possible.
Research of the second author supported by the
American Institute of Mathematics and the NSF Focused Research
Group grant DMS 0244660. The third author was supported by an Overseas Research Scholarship
and a University of Bristol Postgraduate Research Scholarship.
The last author was supported by an EPSRC Advanced Research Fellowship.}

\abstract We prove that a Dirichlet series with a functional
equation and Euler product of a particular form can only arise
from a holomorphic cusp form on the Hecke congruence group
$\Gamma_0(13)$.  The proof does not assume a functional equation
for the twists of the Dirichlet series.  The main new ingredient
is a generalization of the familiar Weil's lemma that played a
prominent role in previous converse theorems.
\endabstract

\maketitle

\section{Introduction and statement of theorem}

An important question in the theory of $L$-functions, is whether a
Dirichlet series with  functional equation and  Euler product of
appropriate type
 can  arise only from some kind of a transform of a related automorphic form.  An affirmative answer to this
question has been given for the simplest types of Dirichlet series
-- those with `degree one' functional equations and arbitrary
conductor,  and degree two functional equations with small
conductors; see the work of Hamburger, Kacorowski-Perelli, and
Hecke~\cite{Ha,H1,KP1,KP2}. In each of these cases, the main ingredient
was the functional equation, the Euler product playing at most a
small role. Conrey and Farmer~\cite{CF} investigated this question in
the setting of Dirichlet series with degree two functional
equations and slightly larger conductors.  For these, it can be
shown that some assumption beyond a functional equation is
absolutely necessary. Weil~\cite{W}, in his converse theorem, imposed
the extra assumption that twists of the given Dirichlet series
also had functional equations. In~\cite{CF}, the more natural condition
that the Dirichlet series has an Euler product -- of the type that
one finds associated to holomorphic modular forms -- is assumed.
They prove that for conductors 5 through 17 (conductors 1 through
4 having been settled by Hecke as mentioned above), with the
possible exception of  13, that all such Dirichlet series are, in
fact, transforms of modular forms.

 In this paper, we introduce a new idea that allows us to fill the gap at 13 in  the
theorem of~\cite{CF}.  The new ingredient (which is in section 5) may
be regarded as a generalization of Weil's lemma, that holomorphic
functions which transform in a certain way under elliptic
transformations of infinite order are identically zero, which
played an important role in~\cite{W, CF}.

Here is a statement of our theorem. Though the notation is
standard, an explanation of it is given later. Also, this paper
almost completely self-contained; some standard arguments are
repeated here for the convenience of the reader.  Below we use the
notation $e(z)=e^{2 \pi i z}$.

\begin{theorem}
\label{theo:1} Suppose
\begin{eqnarray*}
f(z):=\sum_{n=1}^\infty a_n e(nz)
\end{eqnarray*}
is holomorphic in $\Im z>0.$ Suppose further that we have a positive even integer
$k$ such that
\begin{eqnarray*}
L_f(s):=\sum_{n=1}^\infty \frac{a_n}{n^s}
\end{eqnarray*}
converges in some half-plane $\Re s>c$ and that
\begin{eqnarray}\label{eqn:fe}
L_f(s)=\left(1-\frac{a_2}{2^s}+\frac{2^{k-1}}{2^{2s}}\right)^{-1}
\left(1-\frac{a_3}{3^s}+\frac{3^{k-1}}{3^{2s}}\right)^{-1}\sum_{(n,6)=1}\frac{a_n}{n^s}.
\end{eqnarray}
 In other words, we are assuming that the sequence $(a_n)$ doesn't grow too fast and that
it is (degree 2) multiplicative with respect to the primes 2 and 3
and weight $k$. Suppose finally that
\begin{eqnarray*}\Lambda(s)=\left(\frac{\sqrt{13}}{2\pi}\right)^s\Gamma(s)L_f(s)
\end{eqnarray*}
is an entire function which is bounded in any fixed vertical strip, and that it satisfies the functional equation
\begin{eqnarray}\label{eqn:fe2}
\Lambda(s)=\epsilon \Lambda(k-s)
\end{eqnarray}
where $\epsilon=\pm 1$.
Then $f$ is a cusp form of weight $k$ and level 13; i.e. $f\in S_k(\Gamma_0(13)).$

\end{theorem}
 \section{Some notation}
For the convenience of the reader we recall some notation, beginning with the notion of the
 ``stroke'' operator. Let $\gamma=\left(\begin{array}{cc}a&b\\c&d\end{array}\right)$
be a real $2\times 2$ matrix with positive determinant.
 Then
\begin{eqnarray*}
f(z)\vert_k \;\gamma=(\det \gamma)^{k/2}(c
z+d)^{-k}f\left(\frac{az+b}{cz+d}\right).
\end{eqnarray*}
Since $k$ is fixed throughout the paper, we will suppress the dependence on $k$ in this stroke notation.
Also, we will assume that all matrices have positive determinants and real entries. It is easy to verify that
\begin{eqnarray*}
f(z)\vert (\gamma_1 \gamma_2)=(f(z)\vert \gamma_1)\vert \gamma_2
\end{eqnarray*}
and that
\begin{eqnarray*}
f(z)\left\vert \left(\begin{array}{cc}ra&rb\\rc&rd\end{array} \right)=
f(z)\left\vert \left(\begin{array}{cc}a&b\\c&d\end{array} \right)\right.\right.
\end{eqnarray*}
for any real number $r\neq 0$.

To prove Theorem \ref{theo:1}, we need to show
$ f(z)\vert\gamma=f(z) $
for all
\begin{equation*}
\gamma\in \Gamma_0(13):=\left\{\left(\begin{array}{cc}a&b\\c&d\end{array}\right):a,b,c,d\in \mathbb{Z}, ad-bc=1,c\equiv 0\bmod
13\right\}
\end{equation*}
and that $f(z)$ vanishes at all of the cusps of $\Gamma_0(13)$; this is what is meant by $f\in S_k(\Gamma_0(13))$.
It is convenient to
work in the group ring $G=\mathbb{C}[\text{GL}_2^+(\mathbb{R})]$ of formal linear combinations
 of matrices with real entries and positive determinants.
We extend the stroke notation linearly so that
\begin{eqnarray*}
f(z)\vert(a_1 \gamma_1+a_2\gamma_2)=a_1f(z)\vert \gamma_1+a_2f(z)\vert \gamma_2
\end{eqnarray*}
for complex numbers $a_1$ and $a_2$ and real matrices $\gamma_1$ and $\gamma_2$ with positive determinants.
Let $\Omega=\Omega_f=\{\omega\in G:f\vert\omega =0\}$. Then $\Omega$ is a right ideal. It is convenient
to work with
congruences modulo $\Omega$: thus we write
\begin{eqnarray*}\omega_1\equiv \omega_2 \bmod \Omega_f
\end{eqnarray*}
to mean that
\begin{eqnarray*}
f(z)\vert\omega_1 = f(z)\vert \omega_2.
\end{eqnarray*}
To simplify the notation we will usually omit the $\bmod \Omega_f$ from what
we write.
So to prove Theorem~1 we need to verify that
$\gamma\equiv 1$ for all $\gamma\in \Gamma_0(13)$.

Since $\Omega_f$ is  a right ideal one can multiply on the right
a given congruence by anything: thus, $\omega_1\equiv \omega_2$ implies
$\omega_1\omega\equiv \omega_2\omega$ for any $\omega\in G$.
 Also  $\left(\begin{array}{cc}ra&rb\\rc&rd\end{array} \right)\equiv
 \left(\begin{array}{cc}a&b\\c&d\end{array} \right).$

It is not difficult to check that $\Gamma_0(13)$ is generated by four matrices:
\begin{eqnarray*}
\Gamma_0(13)=\left\langle P=\left(\begin{array}{cc}1&1\\0&1\end{array}\right),
W=\left(\begin{array}{cc}1&0\\13&1\end{array}\right),
g_2=\left(\begin{array}{cc}2&-1\\13&-6\end{array}\right),
g_3=\left(\begin{array}{cc}3&-1\\13&-4\end{array}\right) \right\rangle
\end{eqnarray*}
So the main step to prove Theorem \ref{theo:1} is to show $P\equiv
W\equiv g_2\equiv g_3\equiv 1$. The vanishing at the cusps of
$\Gamma_0(13)$ will follow easily, as described near the beginning
of the next section,

\section{Invariance under $P$, $W$, and $g_2$}

Now $P\equiv 1$ asserts exactly the same thing as $f(z+1)=f(z)$,
which follows from the definition of $f(z)$ as a Fourier series.

By Hecke's work, the functional equation (\ref{eqn:fe2}) is
equivalent to $H\equiv \epsilon $ where
\begin{eqnarray*}
H:=\left(\begin{array}{cc}0&-1\\13&0\end{array}\right).\end{eqnarray*}
Since
\begin{eqnarray*}
H\cdot P^{-1}\cdot H =\left(\begin{array}{cc}-13&0\\-169&-13\end{array}\right)\equiv W
\end{eqnarray*}
and $\epsilon^2=1$, we have $W\equiv 1$. That takes care of two of
the four generators of~$\Gamma_0(13)$.

Now we can address the vanishing of $f(z)$ at the cusps.  By the
Fourier series, $f(z)$ vanishes at the cusp~$\infty$.  Since
$f(z)|H=\epsilon f(z)$, and $H$ switches $0$ and~$\infty$, we see
that $f(z)$ also vanishes at~$0$.  But $0$ and $\infty$ are the
only cusps of~$\Gamma_0(13)$, so from the Fourier expansion and
the matrix~$H$, if $f(z)$ is invariant under~$\Gamma_0(13)$ then
$f(z)$ must actually be a cusp form on~$\Gamma_0(13)$.

To prove $g_2\equiv 1$ we need the multiplicativity of $a_n$ at
the prime~$2$.  The following lemma is well-known.

\begin{lemma} We have
\begin{eqnarray*}
\sum_{n=1}^\infty
\frac{a_n}{n^s}=\left(1-\frac{a_p}{p^s}+\frac{p^{k-1}}{p^{2s}}\right)^{-1}\sum_{(n,p)=1}\frac{a_n}{n^s},
\end{eqnarray*}
if and only if
\begin{eqnarray}\label{eqn:A}
\left(\begin{array}{cc}p&0\\0&1\end{array}\right)+\sum_{a=0}^{p-1}\left(\begin{array}{cc}1&a\\0&p\end{array}\right)
\equiv a_p p^{1-{k/2}} .
\end{eqnarray}
\end{lemma}
\begin{proof}
It is convenient to adopt the convention that $a_x=0$ if $x$ is not a positive integer. Equating the coefficient of
$(pn)^{-s}$ on both sides of the equation
\begin{eqnarray*}
\left(1-\frac{a_p}{p^s}+\frac{p^{k-1}}{p^{2s}}\right)\sum_{n=1}^\infty \frac{a_n}{n^s}=\sum_{(n,p)=1}\frac{a_n}{n^s},
\end{eqnarray*}
we have
\begin{eqnarray} \label{eqn:top}
a_{pn}-a_p a_{n}+p^{k-1}a_{n/p}=0.
\end{eqnarray}
A brief calculation shows that
\begin{eqnarray*}
f(z)\left\vert \left(\begin{array}{cc}p&0\\0&1\end{array}\right)+\sum_{a=0}^{p-1}\left(\begin{array}{cc}1&a\\0&p\end{array}\right)
\right. = p^{k/2}f(pz)+p^{1-k/2}\sum_{n=1}^\infty a_{np} e(nz).
\end{eqnarray*}
Thus, equating the coefficient of $e(nz)$ on both sides of (\ref{eqn:A}), we find that
\begin{eqnarray*}
p^{k/2}a_{n/p}+p^{1-k/2}a_{np}=a_pa_{n}p^{1-k/2}
\end{eqnarray*}
which is equivalent to (\ref{eqn:top}).
\end{proof}
Thus, hypothesis (\ref{eqn:fe}) is equivalent to
\begin{eqnarray}\label{eqn:T2}
\left(\begin{array}{cc}2&0\\0&1\end{array}\right)
+\left(\begin{array}{cc}1&0\\0&2\end{array}\right)+\left(\begin{array}{cc}1&1\\0&2\end{array}\right)\equiv
2^{1-k/2}a_2
\end{eqnarray}
and
\begin{eqnarray}\label{eqn:T3}
\left(\begin{array}{cc}3&0\\0&1\end{array}\right)
+\left(\begin{array}{cc}1&0\\0&3\end{array}\right)+\left(\begin{array}{cc}1&1\\0&3\end{array}\right)
+\left(\begin{array}{cc}1&2\\0&3\end{array}\right)\equiv
3^{1-k/2}a_3.
\end{eqnarray}
We multiply each of these equivalences on the left and right by
$H$. (We can multiply on the left by $H$ because $H\equiv \pm 1$).
Using
$H\cdot\left(\begin{array}{cc}a&b\\c&d\end{array}\right)\cdot
H\equiv \left(\begin{array}{cc}d&-c/13\\-13 b&
a\end{array}\right)$ we find that
\begin{eqnarray}\label{eqn:HT2}
\left(\begin{array}{cc}1&0\\0&2\end{array}\right)
+\left(\begin{array}{cc}2&0\\0&1\end{array}\right)+\left(\begin{array}{cc}2&0\\-13&1\end{array}\right)\equiv
2^{1-k/2}a_2
\end{eqnarray}
and
\begin{eqnarray}\label{eqn:HT3}
\left(\begin{array}{cc}1&0\\0&3\end{array}\right)
+\left(\begin{array}{cc}3&0\\0&1\end{array}\right)+\left(\begin{array}{cc}3&0\\-13&1\end{array}\right)
+\left(\begin{array}{cc}3&0\\-26&1\end{array}\right)\equiv
3^{1-k/2}a_3.
\end{eqnarray}
We subtract (\ref{eqn:T2}) from (\ref{eqn:HT2}) to obtain
\begin{eqnarray*}
\left(\begin{array}{cc}2&0\\-13&1\end{array}\right)\equiv
\left(\begin{array}{cc}1&1\\0&2\end{array}\right),
\end{eqnarray*}
from which we deduce that
\begin{eqnarray} \label{eqn:g2}
g_2=W\cdot
\left(\begin{array}{cc}2&0\\-13&1\end{array}\right)\cdot
\left(\begin{array}{cc}1&1\\0&2\end{array}\right)^{-1} \equiv
\left(\begin{array}{cc}2&0\\-13&1\end{array}\right)\cdot
\left(\begin{array}{cc}1&1\\0&2\end{array}\right)^{-1} \equiv 1.
\end{eqnarray}
To complete the proof of Theorem \ref{theo:1} we need only show
that $g_3\equiv 1$.

\section{Three expressions for $f(z)\vert(1-g_3)$}
Invariance under $g_3$ is more difficult, and requires an analytic
argument. We wish to show $f(z)\vert(1-g_3)=0$, so first we
develop some identities for $f(z)\vert(1-g_3)$.

We subtract (\ref{eqn:T3}) from (\ref{eqn:HT3}) to obtain
\begin{eqnarray}\label{eqn:R3}
\left(\begin{array}{cc}1&1\\0&3\end{array}\right)
+\left(\begin{array}{cc}1&2\\0&3\end{array}\right)\equiv \left(\begin{array}{cc}3&0\\-13&1\end{array}\right)
+\left(\begin{array}{cc}3&0\\-26&1\end{array}\right).
\end{eqnarray}
In this expression, we replace $\left(\begin{array}{cc}3&0\\-26&1\end{array}\right)$ by the equivalent matrix
\begin{eqnarray*}
W\cdot \left(\begin{array}{cc}3&0\\-26&1\end{array}\right)=\left(\begin{array}{cc}3&0\\13&1\end{array}\right);
\end{eqnarray*}
 we replace $\left(\begin{array}{cc}1&2\\0&3\end{array}\right)$ by the equivalent matrix
\begin{eqnarray*}
\epsilon H\cdot P^{-1}\cdot
\left(\begin{array}{cc}1&2\\0&3\end{array}\right)= \epsilon
\left(\begin{array}{cc}0&-3\\13&-13\end{array}\right);
\end{eqnarray*}
and  we replace $\left(\begin{array}{cc}3&0\\-13&1\end{array}\right)$ by the equivalent matrix
\begin{eqnarray*}
\epsilon H\cdot\left(\begin{array}{cc}3&0\\-13&1\end{array}\right)=
\epsilon \left(\begin{array}{cc}13&-1\\39&0\end{array}\right);
\end{eqnarray*}
Thus, (\ref{eqn:R3}) can be rewritten as
\begin{eqnarray*}
\left(\begin{array}{cc}1&1\\0&3\end{array}\right) +\epsilon
\left(\begin{array}{cc}0&-3\\13&-13\end{array}\right)-
\epsilon\left(\begin{array}{cc}13&-1\\39&0\end{array}\right)
-\left(\begin{array}{cc}3&0\\13&1\end{array}\right)\equiv 0.
\end{eqnarray*}
Now, multiply on the right by the inverse of the first matrix to
obtain
\begin{eqnarray}\label{eqn:S3}
1 +\epsilon
\left(\begin{array}{cc}0&-3\\39&-26\end{array}\right)-\epsilon
\left(\begin{array}{cc}39&-14\\117&-39\end{array}\right)
-g_3\equiv 0.
\end{eqnarray}
This expression factors as
\begin{eqnarray} \label{eqn:delta1}
\left(1-g_3\right)\cdot \left(1-\epsilon
\left(\begin{array}{cc}39&-14\\117&-39\end{array}\right)\right)\equiv
0
\end{eqnarray}
This expression is the first of three similar factorizations we will find involving $g_3$.

To obtain the second such expression, we first show that
\begin{eqnarray} \label{eqn:H4}
H\cdot\left( \left(\begin{array}{cc}1&1\\0&4\end{array}\right)+\left(\begin{array}{cc}1&3\\0&4\end{array}\right)\right)\cdot H
\equiv \left(\begin{array}{cc}1&1\\0&4\end{array}\right)+\left(\begin{array}{cc}1&3\\0&4\end{array}\right).
\end{eqnarray}
We derive this expression by first squaring (\ref{eqn:T2}) to obtain
\begin{eqnarray*}
&&2 +P +
\left(\begin{array}{cc}4&0\\0&1\end{array}\right)+\left(\begin{array}{cc}1&0\\0&4\end{array}\right)+
\left(\begin{array}{cc}1&1\\0&4\end{array}\right)+
\left(\begin{array}{cc}2&1\\0&2\end{array}\right)+\left(\begin{array}{cc}1&2\\0&4\end{array}\right)+
\left(\begin{array}{cc}1&3\\0&4\end{array}\right)\\
&&\qquad \equiv 2^{2-k}a_2^2.
\end{eqnarray*}
We can replace the terms
$\left(\begin{array}{cc}2&1\\0&2\end{array}\right)$ and
$\left(\begin{array}{cc}1&2\\0&4\end{array}\right)$ here by using
(\ref{eqn:T2}) twice: once multiplied on the right by
$\left(\begin{array}{cc}1&0\\0&2\end{array}\right)$ and once
multiplied on the right by
$\left(\begin{array}{cc}2&0\\0&1\end{array}\right)$. In this way
we obtain, after some rearrangement,
\begin{eqnarray*}
\left(\begin{array}{cc}1&1\\0&4\end{array}\right)+\left(\begin{array}{cc}1&3\\0&4\end{array}\right)
\equiv
2^{2-k}a_2^2-P-2^{1-k/2}a_2\left(\begin{array}{cc}2&0\\0&1\end{array}\right)-
2^{1-k/2}a_2\left(\begin{array}{cc}1&0\\0&2\end{array}\right).
\end{eqnarray*}
The right-hand-side is unchanged when multiplied on the left and right by $H$ which verifies (\ref{eqn:H4}).

Now (\ref{eqn:H4}) can be rewritten as
\begin{eqnarray}  \label{eqn:H5}
 \left(\begin{array}{cc}1&1\\0&4\end{array}\right)+\left(\begin{array}{cc}1&3\\0&4\end{array}\right)
-\left(\begin{array}{cc}4&0\\-13&1\end{array}\right)-\left(\begin{array}{cc}4&0\\-39&1\end{array}\right)
\equiv 0.
\end{eqnarray}

In this expression, we replace $\left(\begin{array}{cc}4&0\\-39&1\end{array}\right)$ by the equivalent matrix
\begin{eqnarray*}
W\cdot \left(\begin{array}{cc}4&0\\-39&1\end{array}\right)=\left(\begin{array}{cc}4&0\\13&1\end{array}\right);
\end{eqnarray*}
 we replace $\left(\begin{array}{cc}1&3\\0&4\end{array}\right)$ by the equivalent matrix
\begin{eqnarray*}
\epsilon H\cdot P^{-1}\cdot
\left(\begin{array}{cc}1&3\\0&4\end{array}\right)= \epsilon
\left(\begin{array}{cc}0&-4\\13&-13\end{array}\right);
\end{eqnarray*}
and  we replace $\left(\begin{array}{cc}4&0\\-13&1\end{array}\right)$ by the equivalent matrix
\begin{eqnarray*}
\epsilon H\cdot\left(\begin{array}{cc}4&0\\-13&1\end{array}\right)=\epsilon \left(\begin{array}{cc}13&-1\\52&0\end{array}\right);
\end{eqnarray*}
thus, (\ref{eqn:H5}) can be rewritten as
\begin{eqnarray}  \label{eqn:H6}
 \left(\begin{array}{cc}1&1\\0&4\end{array}\right)+\epsilon \left(\begin{array}{cc}0&4\\-13&13\end{array}\right)
-\epsilon \left(\begin{array}{cc}13&-1\\52&0\end{array}\right)-\left(\begin{array}{cc}4&0\\13&1\end{array}\right)
\equiv 0.
\end{eqnarray}
Now we multiply on the right by $\left(\begin{array}{cc}1&0\\-13&4\end{array}\right)$; this yields
\begin{eqnarray}  \label{eqn:H7}
 g_3+\epsilon \left(\begin{array}{cc}-26&8\\-91&26\end{array}\right)
-\epsilon \left(\begin{array}{cc}13&-2\\26&0\end{array}\right)-1
\equiv 0.
\end{eqnarray}
This expression factors as
\begin{eqnarray} \label{eqn:delta3}
-\left(1-g_3\right)\cdot \left(1-\epsilon
\left(\begin{array}{cc}-26&8\\-91&26\end{array}\right)\right)\equiv
0
\end{eqnarray}
and gives our second relation of this sort.

The third uses the relation
$g_2g_3^{-1}\equiv g_3 g_2^{-1}$.
This is true because
$g_3^{-1}g_2=\left(\begin{array}{cc}5&-2\\13&-5\end{array}\right)$
has order~2, so $g_2 g_3^{-1}\equiv (g_2 g_3^{-1})^{-1} = g_3 g_2^{-1}$.
Using this relation and $g_2\equiv 1$ we have
\begin{eqnarray*}
1-g_3\equiv 1-g_2 g_3^{-1}g_2\equiv 1-g_3^{-1}g_2\equiv
g_2-g_3^{-1}g_2=-(1-g_3)g_3^{-1}g_2
\end{eqnarray*}
so
\begin{eqnarray}\label{eqn:delta2}
(1-g_3)(1+g_3^{-1}g_2)\equiv 0.
\end{eqnarray}

 \section{Invariance under $g_3$}
In this section, we give an analytic argument to show that
$f$ is invariant under $g_3$.

Let $g(z)=f(z)\vert(1-g_3)$ and let
\begin{eqnarray*}
\delta_1=\left(\begin{array}{cc}\sqrt{13} & \frac{-14}{3\sqrt{13}}\\3 \sqrt{13}& -\sqrt{13}
\end{array}\right),
\qquad
\delta_2=\left(\begin{array}{cc}5&-2\\13&-5\end{array}\right)
\quad {\rm and} \quad \delta_3=\left(\begin{array}{cc}-\sqrt{13} &
\frac{4}{\sqrt{13}}\\\frac{-7 \sqrt{13}}{2}&
\sqrt{13}\end{array}\right).
\end{eqnarray*}
Then by (\ref{eqn:delta1}), (\ref{eqn:delta2}), and (\ref{eqn:delta3}) we have shown that
\begin{equation}\label{eqn:threedeltas}
g(z)\vert \delta_1=\epsilon g(z)\qquad g(z)\vert\delta_2=-g(z) \qquad g(z)\vert\delta_3=\epsilon g(z).
\end{equation}

We will now prove that these relations, \emph{and the fact that
$g_3$ is elliptic} imply that $g(z)$ is~$0$.  The key fact we will
use about $\delta_1$, $\delta_2$ and $\delta_3$ is that
$h_2:=\delta_2\delta_1$ and $h_3:=\delta_3\delta_1$ are irrational
powers of each other.  Therefore $h_2$ and $h_3$ generate a
nondiscrete subgroup of $SL(2,\mathbb R)$, and $g(z)$ is invariant
under stroking by the elements of that group. It would be nice if
this implied that~$g(z)$ is identically zero.  Unfortunately, this
is not quite true, as the following example shows: if
$p(z)=z^{-k/2}$ then $p(z) \vert \left(\begin{array}{cc}X& \\
&1/X\end{array}\right) =p(z)$, for any~$X\in\mathbb R$. This is
essentially the only counterexample, as described in the following
lemma.

\begin{lemma}
If $p(z)$ is an analytic function and
\begin{eqnarray}
p(z)\bigg|   \begin{pmatrix}
        X & 0\\
        0 & 1/X
    \end{pmatrix}
=p(z)
\label{thelemma}
\end{eqnarray}
for all $X$  in a dense subset of $\mathbb{R}_+$,
then $p(z)=C z^{-k/2}$ for some constant~$C$.
\end{lemma}
\begin{proof}
Let $\ell(z)= z^{k/2} p(z)$.  A calculation verifies that
$\ell(z)=\ell(X^2 z)$, so $\ell(z)$ is constant.
\end{proof}

We now put $h_2$ and $h_3$ in a form where we can apply the lemma.
One can check that $h_2$ and $h_3$ commute, so they are simultaneously
diagonalizable.  We have
\begin{eqnarray*}
h_2:=\delta_2\delta_1=\left(\begin{array}{cc}
-\sqrt{13}&\frac{8}{3\sqrt{13}}\\-2\sqrt{13}&\frac{\sqrt{13}}{3}\end{array}\right)
=A \left(\begin{array}
{cc}\frac{-2-\sqrt{13}}{3}&0\\0&\frac{2-\sqrt{13}}{3}\end{array}\right)
A^{-1}\end{eqnarray*}
 and
\begin{eqnarray*}
h_3:=\delta_3\delta_1=\left(\begin{array}{cc}-1 & \frac 23\\
\frac{-13}{2}& \frac{10}{3}\end{array}\right)=A
\left(\begin{array}
{cc}\frac{7-\sqrt{13}}{6}&0\\0&\frac{7+\sqrt{13}}{6}\end{array}\right)
A^{-1}
\end{eqnarray*}
where
\begin{eqnarray*}
A=\left(\begin{array}{cc}  \frac{13+\sqrt{13}}{39}&\frac{13-\sqrt{13}}{39}\\1&1
\end{array}\right).
\end{eqnarray*}
Thus
\begin{eqnarray*}
h_2^m h_3^n =(-1)^mA\left(\begin{array}
{cc}\left(\frac{2+\sqrt{13}}{3}\right)^m&0\\0&\left(\frac{-2+\sqrt{13}}{3}\right)^m\end{array}\right)\left(\begin{array}
{cc}\left(\frac{7-\sqrt{13}}{6}\right)^n&0\\0&\left(\frac{7+\sqrt{13}}{6}\right)^n\end{array}\right)A^{-1}.
\end{eqnarray*}

Let $\lambda=-0.91177\dots$ be the real number such that
\begin{eqnarray*}
 \left(\frac{2+\sqrt{13}}{3}\right)^\lambda=\frac{7-\sqrt{13}}{6}
 \end{eqnarray*}
and let $\displaystyle{Y=\frac{2+\sqrt{13}}{3} }$.
 Then
\begin{eqnarray*}
h_2^m h_3^n =(-1)^mA\left(\begin{array}
{cc}Y^{m+n\lambda}&0\\0&1/Y^{m+n\lambda}
\end{array}\right)
A^{-1}.
\end{eqnarray*}
We have $g(z)|h_2^mh_3^n=(-\epsilon)^m g(z)$ since the number of
$(-1)$'s that we get is the same as the number of times that
$\delta_2$ appears in $h_2^mh_3^n$ and the number of $\epsilon$'s
is the combined number of times that $\delta_1$ and $\delta_3$
appear, which is $m+2n$.

%
%

Replacing $m$ by $2m$ we have $g(z)|h_2^{2m}h_3^n=g(z)$ for all integers~$m, n$.
Thus
$p(z):=g(z)\vert A$ satisfies
\begin{equation}
p(z)|\left(\begin{array}
{cc}X&0\\0&1/X
\end{array}\right)
=p(z)
\end{equation}
for all $X$ of the form $Y^{2m+n\lambda}$ for some integers~$m, n$.
Since $Y$ and $\lambda$ are irrational, the set of such $X$ is dense
in $\mathbb R_+$
and we apply Lemma~\ref{thelemma}
to conclude that~$p(z)=C z^{-k/2}$ for some constant~$C$.  We must
show that~$C=0$.

At this point we must use more information about the
function~$g(z)$.  Indeed, if we let $\tilde g(z) := C z^{-k/2}
|A^{-1}$ then a direct calculation shows for any $C$ that
\begin{equation}
\tilde g(z)\vert \delta_1=(-1)^{-k/2} \tilde g(z)\qquad
\tilde g(z)\vert\delta_2=(-1)^{-k/2} \tilde g(z) \qquad
\tilde g(z)\vert\delta_3= (-1)^{-3k/2} \tilde g(z).
\end{equation}
Thus, if $k\equiv 2 \bmod 4$ and $\epsilon= -1$ then $\tilde g(z)$
satisfies~\eqref{eqn:threedeltas}, so~\eqref{eqn:threedeltas} is
not sufficient by itself to imply that~$g(z)$ is zero. We must use
the fact that $g(z)=f(z)\vert (1-g_3)$ and $g_3$ is elliptic.

Since $g_3^3=I$ we have $(1-g_3)(1+g_3+g_3^2)\equiv 0$. In
particular, $f(z)\vert (1-g_3)(1+g_3+g_3^2)=0$ so
$g(z)|(1+g_3+g_3^2)=0$.  Combining this with the fact that $g(z)=C
z^{-k/2}|A^{-1}$ we obtain
\begin{align*}
0=\mathstrut& C z^{-k/2} \vert \left( 1 + A^{-1}g_3A+A^{-1}g_3^2 A\right) \cr
=\mathstrut & C z^{-k/2}
+
 C z^{-k/2} \Big\vert
\left(
\begin{array}{cc}
 -\frac{1}{2} & \frac{5-2 \sqrt{13}}{6} \\
 \frac{5+2 \sqrt{13}}{6} & -\frac{1}{2}
\end{array}
\right)
+ C z^{-k/2} \Big\vert
\left(
\begin{array}{cc}
 \frac{1}{2} & \frac{5-2 \sqrt{13}}{6} \\
 \frac{5+2 \sqrt{13}}{6} & \frac{1}{2}
\end{array}
\right) \cr
=\mathstrut&
C\bigg( z^{-k/2}
+ 6^k (-3 z-2 \sqrt{13}+5)^{-k/2} ((5+2 \sqrt{13}) z-3)^{-k/2} \cr
&\phantom{XXxXX}
+ 6^k (3 z-2 \sqrt{13}+5)^{-k/2} ((5+2 \sqrt{13}) z+3)^{-k/2}
\bigg).
\end{align*}
The final expression above must be identically~$0$.
Since we assumed $k$ was a positive integer,
if~$C\neq 0$
the final expression above blows
up as $z\to 0$.  Thus $C=0$, so $g(z)=0$, so $f(z)|(1-g_3)=0$,
giving invariance under the final generator of~$\Gamma_0(13)$ and
completing the proof of Theorem~\ref{theo:1}.

It is curious that if $k=-2$ then the final displayed equation above
actually is identically zero.  So the assumption that the
weight $k$ is positive is necessary for the final step of our proof.

\end{document}